\input amssym.def
\input amssym.tex

\hsize=13.5truecm
\baselineskip=16truept plus.3pt minus .3pt

\font\secbf=cmb10 scaled 1200
\font\eightrm=cmr8
\font\sixrm=cmr6

\font\eighti=cmmi8

\font\sixi=cmmi6
\skewchar\eighti='177 \skewchar\sixi='177

\font\eightsy=cmsy8
\font\sixsy=cmsy6
\skewchar\eightsy='60 \skewchar\sixsy='60

\font\eightit=cmti8

\font\eightbf=cmbx8
\font\sixbf=cmbx6

\let\sc=\tensc

\font\eightsc=cmcsc10 scaled 800
\font\secbf=cmb10 scaled 1200
\font\subsecfont=cmb10 scaled \magstephalf
\font\amb=cmmib10

\font\ambi=cmmib10 scaled 700

\newfam\mbfam 

\textfont\mbfam\amb \scriptfont\mbfam\ambi


\def\aa{\def\rm{\fam0\eightrm}%
  \textfont0=\eightrm \scriptfont0=\sixrm \scriptscriptfont0=\fiverm
  \textfont1=\eighti \scriptfont1=\sixi \scriptscriptfont1=\fivei
  \textfont2=\eightsy \scriptfont2=\sixsy \scriptscriptfont2=\fivesy
  \textfont3=\tenex \scriptfont3=\tenex \scriptscriptfont3=\tenex
  \def\sc{\eightsc}
  \def\it{\fam\itfam\eightit}%
  \textfont\itfam=\eightit
  \def\bf{\fam\bffam\eightbf}%
  \textfont\bffam=\eightbf \scriptfont\bffam=\sixbf
   \scriptscriptfont\bffam=\fivebf
  \normalbaselineskip=9.7pt
  \setbox\strutbox=\hbox{\vrule height7pt depth2.6pt width0pt}%
  \normalbaselines\rm}

\def\Proof{\vskip12pt\noindent{\bf Proof.} }
\def\Definition#1{\vskip12pt\noindent{\bf Definition #1}}

\def\Remark#1{\vskip12pt\noindent{\bf Remark #1}}

\def\m@th{\mathsurround=0pt}

\def\cc#1{\hbox to .89\hsize{$\displaystyle\hfil{#1}\hfil$}\cr}
\def\lc#1{\hbox to .89\hsize{$\displaystyle{#1}\hfill$}\cr}
\def\rc#1{\hbox to .89\hsize{$\displaystyle\hfill{#1}$}\cr}

\def\eqal#1{\null\,\vcenter{\openup\jot\m@th
  \ialign{\strut\hfil$\displaystyle{##}$&&$\displaystyle{{}##}$\hfil
      \crcr#1\crcr}}\,}

\def\section#1{\vskip 22pt plus6pt minus2pt\penalty-400
        {{\secbf
        \noindent#1\rightskip=0pt plus 1fill\par}}
        \par\vskip 12pt plus5pt minus 2pt
        \penalty 1000}

\def\subsection#1{\vskip 20pt plus6pt minus2pt\penalty-400
        {{\subsecfont
        \noindent#1\rightskip=0pt plus 1fill\par}}
        \par\vskip 8pt plus5pt minus 2pt
        \penalty 1000}

\def\subsubsection#1{\vskip 18pt plus6pt minus2pt\penalty-400
        {{\subsecfont
        \noindent#1}}
        \par\vskip 7pt plus5pt minus 2pt
        \penalty 1000}

\def\centerlast#1{\begingroup \leftskip=0pt plus 1fil\rightskip=0pt plus -1fil
\parfillskip=0pt plus 2fil
\parindent 0pt
\par #1\par\endgroup}
\def\\{\hfill\break}

\def\kwadrat{\null\ \hfill\null\ \hfill$\square$}
\def\mida#1{{{\null\kern-4.2pt\left\bracevert\vbox to 6pt{}\!\hbox{$#1$}\!\right\bracevert\!\!}}}
\def\midy#1{{{\null\kern-4.2pt\left\bracevert\!\!\hbox{$\scriptstyle{#1}$}\!\!\right\bracevert\!\!}}}

\def\diagin{\diagup\hskip-11.5pt\intop}
\def\diagint{{\raise1.5pt\hbox{$\scriptscriptstyle\diagup$}\hskip-8.7pt\intop}}
\def\diagintop{\mathop{\mathchoice
{{\diagin}}%
{{\diagint}}%
{{\diagint}}%
{{\diagint}}%
}\limits}

\def\divv{{{\rm div}\,}}

\def\today{${\scriptscriptstyle\number\day-\number\month-\number\year}$}
\footline={{\hfil\rm\the\pageno\hfil${\scriptscriptstyle\rm\jobname}$\ \ \today}}

\def\ifnextchar#1#2#3{\bgroup
  \def\reserveda{\ifx\reservedc #1 \aftergroup\firstoftwo
    \else \aftergroup\secondoftwo\fi\egroup{#2}{#3}}%
  \futurelet\reservedc\ifnch
  } 
\def\ifnch{\ifx \reservedc \sptoken \expandafter\xifnch
      \else \expandafter\reserveda
      \fi} 
\def\firstoftwo#1#2{#1}
\def\secondoftwo#1#2{#2} 
\def\tempswafalse{\let\iftempswa\iffalse}
\def\tempswatrue{\let\iftempswa\iftrue} 

\def\cite{\ifnextchar [{\tempswatrue\citea}{\tempswafalse\citeb}}
\def\citea[#1]#2{[#2, #1]}
\def\citeb#1{[#1]}\phantom{]}

\def\\{\hfil\break}

\def\N{{\Bbb N}}
\def\R{{\Bbb R}}
\def\T{{\Bbb T}}

\def\divv{{\rm div}\,}

\centerlast{\secbf Stability of two-dimensional Navier-Stokes motions in the 
periodic case}

\vskip1cm
\centerline{\bf E. Zadrzy\'nska$^1$ and W. M. Zaj\c aczkowski$^2$}
\vskip1cm
\item{$^1$} Faculty of Mathematics and Information Sciences, Warsaw University of Technology, Koszykowa 75, 00-662 Warsaw, Poland.\\ 
e-mail:emzad@mini.pw.edu.pl
\item{$^2$} Institute of Mathematics Polish Academy of Sciences,\\ 
\'Sniadeckich 8, 00-656 Warsaw, Poland,\\
e-mail:wz@impan.gov.pl;\\
Institute of Mathematics and Cryptology, Cybernetics Faculty,\\ 
Military University of Technology, Kaliskiego 2, 00-908 Warsaw, Poland\vskip1cm

\noindent
{\bf Abstract.} 
We consider the motion described by the Navier-Stokes equations in a box with 
periodic boundary conditions. First we prove the existence of global strong 
two-dimensional solutions. Next we show the existence of global strong 
three-dimensional solutions under the assumption that the initial data and 
the external force are sufficiently close to the initial data and the external 
force of the two-dimensional problem in appropriate spaces. The second 
result can be treated as stability of strong two-dimensional solutions in 
the set of suitably strong three-dimensional motions.

MSC 2010: 35Q30, 76D05, 76N10, 35B35, 76D03

Key words: incompressible Navier-Stokes equations, stability of 
two-dimensional solutions, global regular solutions

\section{1. Introduction}

The aim of this paper is to prove stability of two-dimensional periodic 
solutions in the set of three-dimensional periodic solutions to the 
Navier-Stokes equations. We consider three-dimensional fluid motions in a box 
$\Omega=[0,L]^3$, $L>0$, described by
$$\eqal{
&v_t+v\cdot\nabla v-\nu\Delta v+\nabla p=f\quad &{\rm in}\ \ 
\Omega\times\R_+,\cr
&\divv v=0\quad &{\rm in}\ \ \Omega\times\R_+,\cr
&v|_{t=0}=v(0)\quad &{\rm in}\ \ \Omega,\cr}
\leqno(1.1)
$$
where $v=(v_1(x,t),v_2(x,t),v_3(x,t))\in\R^3$ is the velocity of the fluid, 
$x=(x_1,x_2,x_3)$ with $x_i\in(0,L)$, $i=1,2,3$, $p=p(x,t)\in\R$ is the 
pressure and $f=(f_1(x,t),f_2(x,t),f_3(x,t))\in\R^3$ is the external force 
field.

\noindent
Finally, $\nu>0$ is the constant viscosity coefficient and the dot in the 
second term of $(1.1)_1$ denotes the scalar product.

By two-dimensional motions we mean solutions $(v,p)$ to (1.1) such that 
$v=v_s=\break(v_{s1}(x_1,x_2,t)$, $v_{s2}(x_1,x_2,t),0)\in\R^2$, 
$p=p_s(x_1,x_2,t)\in\R$ and 
$f=f_s=(f_{s1}(x_1,x_2,t),\break f_{s2}(x_1,x_2,t),0)\in\R^2$.

\noindent
The main result of this paper is the following. Assume that $f-f_s$ and 
$v(0)-v_s(0)$ are sufficiently small in some norms. Then we show that $v-v_s$ 
and $p-p_s$ are small in appropriate norms for all times. Observe that we are 
talking about global solutions.

More precisely, two-dimensional periodic solutions satisfy
$$\eqal{
&v_{st}+v_s\cdot\nabla v_s-\nu\Delta v_s+\nabla p_s=f_s\quad &{\rm in}\ \ 
\Omega\times\R_+,\cr
&\divv v_s=0\quad &{\rm in}\ \ \Omega\times\R_+,\cr
&v_s|_{t=0}=v_s(0)\quad &{\rm in}\ \ \Omega,\cr}
\leqno(1.2)
$$
where no quantities in (1.2) depend on $x_3$. To show stability, 
we introduce the quantities
$$
u=v-v_s,\quad q=p-p_s
$$
which are periodic solutions to the problem
$$\eqal{
&u_t+u\cdot\nabla u-\nu\Delta u+\nabla q=-v_s\cdot\nabla u-u\cdot\nabla v_s+g
\quad &{\rm in}\ \ \Omega\times\R_+,\cr
&\divv u=0\quad &{\rm in}\ \ \Omega\times\R_+,\cr
&u|_{t=0}=u(0)\quad &{\rm in}\ \ \Omega,\cr}
\leqno(1.3)
$$
where $g=f-f_s$.

\noindent
Our aim is to show the smallness of $u(t)$ for all $t\in\R_+$ if $u(0)$ and 
$g$ are sufficiently small.

\noindent
To derive necessary estimates we use the energy method. Hence the Poincar\'e 
inequality is needed. But it does not hold for solutions to problems (1.2) 
and (1.3). Therefore we introduce the quantities
$$\eqal{
&\bar v_s=v_s-\diagintop_\Omega v_sdx,\quad 
&\bar u=u-\diagintop_\Omega udx,\quad &\bar p_s=p_s-\diagintop_\Omega p_sdx,\cr
&\bar q=q-\diagintop_\Omega qdx,\quad &\bar f_s=f_s-\diagintop_\Omega f_sdx,\quad
&\bar g=g-\diagintop_\Omega gdx,\cr}
\leqno(1.4)
$$
where the integral mean is defined by
$$
\diagintop_\Omega\omega dx={1\over|\Omega|}\intop_\Omega\omega dx
$$
and $|\Omega|=L^3$. For the quantities (1.4) problems (1.2) and (1.3) take 
the form
$$\eqal{
&\bar v_{st}+v_s\cdot\nabla\bar v_s-\nu\Delta\bar v_s+\nabla\bar p_s=\bar f_s
\quad &{\rm in}\ \ \Omega\times\R_+,\cr
&\divv\bar v_s=0\quad &{\rm in}\ \ \Omega\times\R_+,\cr
&\bar v_s|_{t=0}=\bar v_s(0)\quad &{\rm in}\ \ \Omega,\cr}
\leqno(1.5)
$$
and
$$\eqal{
&\bar u_t+u\cdot\nabla\bar u-\nu\Delta\bar u+\nabla\bar q=-v_s\cdot\nabla
\bar u-u\cdot\nabla\bar v_s+\bar g\quad &{\rm in}\ \ \Omega\times\R_+,\cr
&\divv\bar u=0\quad &{\rm in}\ \ \Omega\times\R_+,\cr
&\bar u|_{t=0}=\bar u(0)\quad &{\rm in}\ \ \Omega.\cr}
\leqno(1.6)
$$
For functions $\bar v_s$, $\bar u$ the Poincar\'e inequality does hold.

\noindent
Since we are looking for periodic solutions to problem (1.1) we introduce the 
notation: $H^m(\Omega)=\{u\in H_{loc}^m(\R^3):\ u(x+Le_i)=u(x),i=1,2,3\}$, 
where $e_i$, $i=1,2,3$ is the canonical basis and
$$
H^{2,1}(\Omega\times I)=\{u=u(x,t):\ u\in L_2(I;H^2(\Omega)),
u_t\in L_2(I;L_2(\Omega))\},
$$
where $I\subset\R$ is an open interval.

Let $\bar v_s(0)\in H^1(\Omega)$, $f_s\in L_{2,loc}(\R_+;H^1(\Omega))$ and 
denote
$$
\bar A_1^2(T)=\sup_{k\in\N_0}\intop_{kT}^{(k+1)T}
\|\bar f_s(t)\|_{H^1(\Omega)}^2dt,
$$
where $T>0$, $\N_0=\N\cup\{0\}$;
$$\eqal{
&\bar A_2^2=\|\bar v_s(0)\|_{H^1(\Omega)}^2,\cr
&\bar A_3^2(T)=c_1(\bar A_1^2(T)+\bar A_2^2)\bar A_2^2+(\bar A_1^2(T)+1)
e^{c_2(\bar A_1^2(T)+\bar A_2^2)},\cr}
$$
where $c_1,c_2>0$ are some constants;
$$
T_*={2\over c_{s1}}\ln2,
$$
where $c_{s1}>0$ is a constant depending on $\nu$ (introduced in Lemma 2.2).

For a given $\bar v_s(0)\in H^1(\Omega)$ define
$$
{\cal M}=\{(T,f_s)\in[T_*,\infty)\times L_{2,loc}(\R_+;H^1(\Omega)):\ 
T>\bar A_3^2(T)\}.
$$

\Definition{1.1.} 
Let $(T,f_s)\in{\cal M}$, $\bar v_s(0)\in H^1(\Omega)$, $\divv\bar v_s(0)=0$. 
A pair of functions $(\bar v_s,\bar p_s)$ is called a strong solution to 
problem (1.5) if $\bar v_s$ is a weak solution of system $(1.5)_{1,2}$ in 
$\Omega\times(kT,(k+1)T)$ with the initial condition 
$\bar v_s|_{t=kT}=\bar v_s(kT)$ for all $k\in\N_0$ and if 
$\bar v_s\in L_\infty(kT,(k+1)T;H^1(\Omega))\cap L_2(kT,(k+1)T;H^2(\Omega))$, 
$\nabla\bar p_s\in L_2(\Omega\times(kT,(k+1)T))$ for all $k\in\N_0$.

\noindent
Analogous definition holds for solutions to problem (1.1).

\proclaim Theorem 1.1. 
Let $\bar v_s(0)\in H^2(\Omega)$, $\divv\bar v_s(0)=0$, 
$f_s\in L_{2,loc}(\R_+;H^1(\Omega))$, $\bar A_1^2(T)<\infty$ for all $T>0$ and 
assume that $(T,f_s)\in{\cal M}$. \\
Then there exists a unique strong solution $(\bar v_s,\bar p_s)$ to problem 
(1.5) such that \\
$\bar v_s\in H^{2,1}(\Omega\times(kT,(k+1)T))\cap C([kT,(k+1)T];H^2(\Omega))\cap
L_2(kT,(k+1)T;H^3(\Omega))$, $\nabla\bar p_s\in L_2(\Omega\times(kT,(k+1)T))$, 
$k\in\N_0$ and 
$$\eqal{
&\|\bar v_s\|_{C([kT,(k+1)T];H^2(\Omega))}^2+
\|\bar v_s\|_{L_2(kT,(k+1)T;H^3(\Omega))}^2
\le c(\bar A_1^2,\bar A_2^2)(\bar A_1^2+\bar A_2^2),\cr}
$$
where $c=c(\bar A_1^2,\bar A_2^2)$ does not depend on $k$.

Notice that the set of admissible functions $f_s$ is large.
Let us give two examples of such functions. First define 
$f_s=a+h_s$, where $a\in\R$ and $h_s\in L_2(\R_+;H^1(\Omega))$. Then 
$\bar f_s=\bar h_s\in L_2(\R_+;H^1(\Omega))$ and if 
$T>c_1\big(\intop_0^\infty\|\bar h_s(t)\|_{H^1(\Omega)}^2dt+\bar A_2^2\big)+
\big(\intop_0^\infty\|\bar h_s(t)\|_{H^1(\Omega)}^2dt+1\big)
e^{c_2(\intop_0^\infty\|\bar h_s(t)\|_{H^1(\Omega)}^2dt+\bar A_2^2)}\equiv 
A_0$ we have $(T,f_s)\in{\cal M}$ for all $T>\max(T_*,A_0)$.

The above example shows that there is no restriction on the magnitude of the 
external force $f_s$ and for $a\not=0$, $f_s$ need not decay in time.

Now, consider another example. Let $h_s\in L_2(\R_+;H^1(\Omega))$ and 
$T>A_0$. Define a periodic function $f_{sT}(x,t)=h_s(x,t-kT)$ for 
$kT\le t\le(k+1)T$, $k\in\N_0$. Then $T>\bar A_3^2(T)$ and 
$(T,f_{sT})\in{\cal M}$ for all $T>\max(T_*,A_0)$.

Theorem 1.1 yields the existence of a solution to problem (1.5) such that 
$\bar v_s\in H^{2,1}(\Omega\times(kT,(k+1)T)\cap C([kT,(k+1)T];
H^2(\Omega))\cap L_2(kT,(k+1)T;H^3(\Omega)))$. However, the assumptions of the 
theorem are too weak to obtain an estimate of 
$\|\bar v_s\|_{H^{2,1}(\Omega\times(kT,(k+1)T))}$ which is independent of $k$. 
To derive such an estimate we need an additional assumption on $\bar f_s$,
formulated in the theorem below.

\proclaim Theorem 1.2. 
Let the assumptions of Theorem 1.1 hold. Moreover, suppose that
$$
\bar A_4^2=\sup_{k\in\N_0}\sup_{kT\le t\le(k+1)T}\bigg|\intop_0^t
\diagintop_\Omega f_s(x,t')dxdt'+\diagintop_\Omega v_s(0)dx\bigg|^2<\infty.
$$
Then there exists a unique strong solution $(\bar v_s,\bar p_s)$ to problem 
(1.5) such that \\
$\bar v_s\in H^{2,1}(\Omega\times(kT,(k+1)T)\cap C([kT,(k+1)T];H^2(\Omega))
\cap L_2(kT,(k+1)T;H^3(\Omega)))$,\break
$\nabla\bar p_s\in L_2(\Omega\times(kT,(k+1)T))$, $k\in\N_0$ and
$$\eqal{
&\|\bar v_{st}\|_{L_2(kT,(k+1)T;L_2(\Omega))}^2+
\|\bar v_s\|_{C([kT,(k+1)T];H^2(\Omega))}^2\cr
&\quad+\|\bar v_s\|_{L_2(kT,(k+1)T;H^3(\Omega))}^2+
\|\nabla\bar p_s\|_{L_2(\Omega\times(kT,(k+1)T))}^2\le c
(\bar A_1^2,\bar A_2^2,\bar A_4^2)(\bar A_1^2+\bar A_2^2),\cr}
$$
where $c=c(\bar A_1^2,\bar A_2^2,\bar A_4^2)$ does not depend on $k$.

Using Theorem 1.1 the following theorem concerning the stability of 
a two-dimensio\-nal solution in the set of three-dimensional solutions can be 
proved. This theorem gives also the existence of a global strong solution 
to problem (1.1)
\eject

\proclaim Theorem 1.3. 
Let the assumptions of Theorem 1.1 hold. Let $v(0)\in H^1(\Omega)$,\\
$\divv v(0)=0$, $f\in L_{2,loc}(\R_+;L_2(\Omega))$ and suppose that
$$\eqal{
\bar G(t)&=\sup_{k\in\N_0}\intop_{kT}^{(k+1)T}\|\bar g(t')\|_{L_2(\Omega)}^2
dt'+\|\bar u(0)\|_{L_2(\Omega)}^2\cr
&\quad+\sup_{k\in\N_0}\intop_{kT}^{(k+1)T}\bigg|\intop_0^t\diagintop_\Omega 
gdxdt'+\diagintop_\Omega u(0)dx\bigg|^2dt\cr
&\quad+\bigg|\intop_0^t\diagintop_\Omega gdxdt'+
\diagintop_\Omega u(0)dx\bigg|^2+\|\bar g(t)\|_{L_2(\Omega)}^2<\infty\quad 
{\sl for\ all}\ \ t\in\R_+.\cr}
$$
There exists a constant $\gamma>0$ such that if
$$
\|u(0)\|_{H^1(\Omega)}^2\le\gamma,\quad \bar G(t)\le\varepsilon\gamma\quad 
{\sl for\ all}\ \ t\in\R_+,\quad {\sl and\ some}\ \ 0<\varepsilon<1,
$$
then there exists a unique strong solution $(v,p)$ to problem (1.1) such that\\
$v\in H^{2,1}(\Omega\times(kT,(k+1)T))$, 
$\nabla p\in L_2(\Omega\times(kT,(k+1)T))$, $k\in\N_0$, and
$$
\|u(t)\|_{H^1(\Omega)}^2\le c\gamma\quad {\sl for\ all}\ \ t\in\R_+,
\leqno(1.7)
$$
where $c>0$ is some constant. Moreover,
$$
\|u\|_{L_2(kT,(k+1)T;H^2(\Omega))}^2\le\bar c\gamma\quad {\sl for\ all}\ \ 
k\in\N_0,
\leqno(1.8)
$$
where $\bar c=\bar c(T)$.

Notice that Theorem 1.3 yields the existence of $v$ in $H^{2,1}$ while the 
stability of $v_s$ in a weaker norm. In the theorem below we formulate 
the stability result for $H^{2,1}$-norm.

\proclaim Theorem 1.4. 
Let the assumptions of Theorems 1.2 and 1.3 be satisfied. Moreover suppose that
$$
\sup_{k\in\N_0}\sup_{t\in[kT,(k+1)T]}\bigg|\intop_0^t\diagintop_\Omega g(t')
dxdt'+\diagintop_\Omega u(0)dx\bigg|^2\le\gamma.
$$
If $\gamma$ is sufficiently small then the solution $(v,p)$ of problem (1.1), 
which exists in virtue of Theorem 1.3, satisfies
$$
\|u\|_{H^{2,1}(\Omega\times(kT,(k+1)T))}^2+
\|\nabla q\|_{L_2(kT,(k+1)T;L_2(\Omega))}^2\le c\gamma,
\leqno(1.9)
$$
where $c=c(T)$, $k\in\N_0$.

The stability problem for Navier-Stokes equations has been developed in 
different directions. There are results concerning the stability of weak or 
regular solutions as well as the stability of two-dimensional solutions or 
other special solutions in the three-dimensional space. Some papers discuss 
the question of stability of stationary solutions in the set of nonstationary 
solutions.

The first results connected with the stability of global regular solutions to 
the nonstationary Navier-Stokes equations were proved by Beirao da Veiga and 
Secchi \cite{2}, followed by Ponce, Racke, Sideris and Titi \cite{13}. 
Paper \cite{2} is concerned with the stability in $L_p$-norm of a strong 
three-dimensional solution of the Navier-Stokes system with zero external 
force in the whole space. In \cite{13}, assuming that the external force is 
zero and a three-dimensional initial function is close to 
a two-dimensional one in $H^1(\R^3)$, the authors showed the existence of 
a global strong solution in $\R^3$ which remains close to a two-dimensional 
strong solution for all times. In \cite{12} Mucha obtained a similar result 
under weaker assumptions about the smallness of the initial velocity 
perturbation.

In the class of weak Leray-Hopf solutions the first stability result was 
obtained by Gallagher \cite{6}. She proved the stability of 
two-dimensional solutions of the Navier-Stokes equations with 
periodic boundary conditions under three-dimensional perturbations both in 
$L_2$ and $H^{1\over2}$ norms.

The stability of nontrivial periodic regular solutions to the Navier-Stokes 
equations was studied by Iftimie \cite{8} and by Mucha \cite{10}. 
The paper \cite{10} is devoted to the case when the external force is 
a potential belonging to $L_{r,loc}(\T^3\times[0,\infty))$ and when the intial 
data belongs to the space $W_r^{2-2/r}(\T^3)\cap L_2(\T^3)$, where $r\ge2$ and 
$\T$ is a torus. Under the assumption that there exists a global solution with 
data of regularity mentioned above and assuming that small perturbations 
of data have the same regularity as above, the author proves that 
perturbations of the velocity and the gradient of the pressure remain small in 
the spaces $W_r^{2,1}(\T^3\times(k,k+1))$ and $L_r(\T^3\times(k,k+1))$, 
$k\in\N$, respectively. Paper \cite{8} contains results concerning the 
stability of two-dimensional regular solutions to the Navier-Stokes system in 
a three-dimensional torus but here the initial data in the three-dimensional 
problem belongs to an anisotropic space of functions having different 
regularity in the first two directions than in the third direction, and the 
external force vanishes. Moreover, Mucha \cite{11} studies the stability of 
regular solutions to the nonstationary Navier-Stokes system in $\R^3$ assuming 
that they tend in $W_r^{2,1}$ spaces $(r\ge2)$ to constant flows.

The papers of Auscher, Dubois and Tchamitchian \cite{1} and of Gallagher, 
Iftimie and Planchon \cite{7} concern the stability of global regular 
solutions to the Navier-Stokes equations in the whole space $\R^3$ with zero 
external force. These authors assume that the norms of the solutions 
considered decay as $t\to\infty$.

It is worth mentioning the paper of Zhou \cite{14}, who proved the asymptotic 
stability of weak solutions $u$ with the property: 
$u\in L_2(0,\infty,BMO)$ to the Navier-Stokes equations in $\R^n$, $n\ge3$, 
with force vanishing as $t\to\infty$.

An interesting result was obtained by Karch and Pilarczyk \cite{9}, who 
concentrate on the stability of Landau solutions to the Navier-Stokes system 
in $\R^3$. Assuming that the external force is a singular distribution they 
prove the asymptotic stability of solution under any $L_2$-perturbation.

Paper \cite{5} of Chemin and Gallagher is devoted to the stability 
of some unique global solution with large data in a very weak sense.

Finally, the stability of Leray-Hopf weak solutions has recently been examined 
by Bardos et al. \cite{3}, where equations with vanishing external force are 
considered. That paper concerns the following three cases: two-dimensional 
flows in infinite cylinders under three-dimensional perturbations which are 
periodic in the vertical direction; helical flows in circular cylinders under 
general three-dimensional perturbations; and axisymmetric flows under general 
three-dimensional perturbations. The theorem concerning the first case 
extends a result obtained by Gallagher \cite{6} for purely periodic 
boundary conditions.

\noindent
Most of the papers discussed above concern to the case with zero external 
force (\cite{1}--\cite{3}, \cite{5}--\cite{8}, \cite{12}, \cite{13}) 
or with force which decays as $t\to\infty$ (\cite{18}). Exceptions are 
\cite{9}--\cite{11}, where very special external forces, which are singular 
distributions in \cite{9} or potentials in \cite{10}--\cite{11}, are 
considered. However, the case of potential forces is easily reduced to the 
case of zero external forces.

The aim of our paper is to prove the stability result for a large class of 
external forces $f_s$ which do not produce solutions decaying as $t\to\infty$. 
Examples of such functions have been given after the formulation 
of Theorem 1.1.

It is essential that our stability results are obtained together with the 
existence of a global strong three-dimensional solution close to 
a two-dimensional one.

The paper is divided into two main parts. In the first we prove existence of 
global strong two-dimensional solutions not vanishing as $t\to\infty$ because 
the external force does not vanish either. To prove existence of such 
solutions we use the step by step method. For this purpose we have to show 
that the data in the time interval $[kT,(k+1)T]$, $k\in\N$, do not increase 
with $k$. For this we also need the time step $T$ to be sufficiently large.

\noindent
In the second part we prove existence of three-dimensional solutions that 
remain close to two-dimensional solutions. For this we need the initial 
velocity and the external force to be sufficiently close in apropriate norms 
to the initial velocity and the external force of the two-dimensional 
problems.

The proofs of this paper are based on the energy method, which is available 
thanks to the periodic boundary conditions. The proofs of global existence 
which follow from the step by step technique are possible thanks to the natural 
decay property of the Navier-Stokes equations. This is mainly used in the 
first part of the paper (Section 3). To prove stability (Section 4) we use 
smallness of data $(v(0)-v_s(0)),(f-f_s)$ and a contradiction argument applied 
to the nonlinear ordinary differential inequality (4.11).

We restrict ourselves to proving estimates,because existence follows easily 
by the Faedo-Galerkin method.

The paper is organized as follows. In Section 2 we introduce notation and give 
some auxiliary results. Section 3 is devoted to the existence of 
a two-dimensional solution. It also contains some useful estimates of the 
solution. In Section 4 we prove the existence of a global strong solution 
to problem (1.1) close to the two-dimensional solution for all time.

\section{2. Notation and auxiliary results}

By $L_p(\Omega)$, $p\in[1,\infty]$, we denote the Lebesgue space of integrable 
functions. By $H^s(\Omega)$, $s\in\N_0=\N\cup\{0\}$, we denote the Sobolev 
space of periodic functions with the finite norm
$$
\|u\|_{H^s}\equiv\|u\|_{H^s(\Omega)}=\sum_{|\alpha|\le s}\bigg(\intop_\Omega
|D_x^\alpha u|^2dx\bigg)^{1/2},
$$
where $D_x^\alpha=\partial_{x_1}^{\alpha_1}\partial_{x_2}^{\alpha_2}
\partial_{x_3}^{\alpha_3}$, 
$|\alpha|=\alpha_1+\alpha_2+\alpha_3$, $\alpha_i\in\N_0$, $i=1,2,3$.

\noindent
To prove Theorems 1.2, 1.4 we need formulas for the means of $v_s$ and $u$. 
Hence, we have

\proclaim Lemma 2.1. 
Assume that $\diagintop_\Omega f_s(t)dx$, $\diagintop_\Omega g(t)dx$ are 
locally integrable on $\R_+$ and $\diagintop_\Omega v_s(0)dx$, 
$\diagintop_\Omega u(0)dx$ are finite. Then, for all $t\in\R_+$,
$$
\diagintop_\Omega v_s(t)dx=\intop_0^t\diagintop_\Omega f_s(t)dx+
\diagintop_\Omega v_s(0)dx,
\leqno(2.1)
$$
$$
\diagintop_\Omega u(t)dx=\intop_0^t\diagintop_\Omega g(t)dx+
\diagintop_\Omega u(0)dx.
\leqno(2.2)
$$

\Proof 
Applying the mean operator to (1.2) and (1.4), integrating by parts and using 
the periodic boundary conditions, we get
$$
{d\over dt}\diagintop_\Omega v_sdx=\diagintop_\Omega f_sdx,
\leqno(2.3)
$$
$$
{d\over dt}\diagintop_\Omega udx=\diagintop_\Omega gdx.
\leqno(2.4)
$$
Integrating (2.3) and (2.4) with respect to time yields (2.1) and (2.2).
\kwadrat

The following lemma follows directly from the Poincar\'e inequality.

\proclaim Lemma 2.2. 
We have
$$
c_{s1}\|\bar v_s\|_{H^1}^2\le\nu\|\nabla\bar v_s\|_{L_2}^2,
\leqno(2.5)
$$
$$
c_1\|\bar u\|_{H^1}^2\le\nu\|\nabla\bar u\|_{L_2}^2,
\leqno(2.6)
$$
where $c_1$, $c_{s1}$ are positive constants.

\section{3. Two-dimensional solutions}

First we need

\proclaim Lemma 3.1. 
Assume that
$$\eqal{
&1.\qquad A_1^2&={1\over c_{s1}}\sup_{k\in\N_0}\intop_{kT}^{(k+1)T}
\|\bar f_s(t')\|_{L_2}^2dt'<\infty,\cr
&2.\qquad A_2^2&={A_1^2\over1-e^{-c_{s1}T}}+\|\bar v_s(0)\|_{L_2}^2<\infty,\cr}
$$
where $T>0$ is fixed and $c_{s1}$ is introduced in (2.5). Then
$$
\|\bar v_s(kT)\|_{L_2}^2\le A_2^2
\leqno(3.1)
$$
and
$$
\|\bar v_s(t)\|_{L_2}^2+c_{s1}\intop_{kT}^t
\|\bar v_s(t')\|_{H^1}^2dt'\le A_1^2+A_2^2\equiv A_3^2
\leqno(3.2)
$$
for all $t\in(kT,(k+1)T]$, $k\in\N_0$.

\Proof 
Multiplying $(1.5)_1$ by $\bar v_s$, integrating over $\Omega$, using the 
periodic boundary conditions and inequality (2.5) yields
$$
{1\over2}{d\over dt}\|\bar v_s\|_{L_2}^2+c_{s1}\|\bar v_s\|_{H^1}^2\le
{c_{s1}\over2}\|\bar v_s\|_{L_2}^2+{1\over2c_{s1}}\|\bar f_s\|_{L_2}^2,
$$
where we also applied the Young inequality to the term with the r.h.s. of 
$(1.5)_1$.

\noindent
Hence, we have
$$
{d\over dt}\|\bar v_s\|_{L_2}^2+c_{s1}\|\bar v_s\|_{H^1}^2\le{1\over c_{s1}}
\|\bar f_{sx}\|_{L_2}^2.
\leqno(3.3)
$$
Continuing, we obtain
$$
{d\over dt}(\|\bar v_s\|_{L_2}^2e^{c_{s1}t})\le{1\over c_{s1}}
\|\bar f_s\|_{L_2}^2e^{c_{s1}t}.
$$
Integrating with respect to time yields
$$
\|\bar v_s(t)\|_{L_2}^2\le{1\over c_{s1}}\intop_{kT}^t
\|\bar f_s(t')\|_{L_2}^2dt'+e^{-c_{s1}(t-kT)}\|\bar v_s(kT)\|_{L_2}^2,
$$
for all $k\in\N_0$, $T>0$ and $t\in(kT,(k+1)T]$. Setting $t=(k+1)T$ we get
$$
\|\bar v_s((k+1)T)\|_{L_2}^2\le{1\over c_{s1}}\intop_{kT}^{(k+1)T}
\|\bar f_s(t')\|_{L_2}^2dt'+e^{-c_{s1}T}\|\bar v_s(kT)\|_{L_2}^2.
$$
By iteration we have
$$
\|\bar v_s(kT)\|_{L_2}^2\le{A_1^2\over1-e^{-c_{s1}T}}+e^{-c_{s1}kT}
\|v_s(0)\|_{L_2}^2\le A_2^2.
$$
Hence, (3.1) is proved. Integrating (3.3) with respect to time from $t=kT$ to 
$t\in(kT,(k+1)T]$, we obtain (3.2).
\kwadrat

\noindent
To obtain an estimate for the second derivatives of $\bar v_s$ we need

\proclaim Lemma 3.2. 
Let the assumptions of Lemma 3.1 hold. Let $\bar v_s(0)\in H^1(\Omega)$, 
$\divv\bar v_s(0)=0$. Suppose that
$$
T\ge{2c_{s2}\over c_{s1}}A_3^2,
\leqno(3.4)
$$
where $c_{s1}$ is the constant from inequality (2.5), $c_{s2}$ is introduced 
in (3.8) below and $A_3^2$ is defined in Lemma 3.1. Denote\\
1. $A_4^2=c_{s1}e^{c_{s1}A_3^2}A_1^2$,\\
2. $A_5^2={A_4^2\over1-e^{-{c_{s1}\over2}T}}+\|\bar v_{sx}(0)\|_{L_2}^2$,\\
3. $A_6^2=A_4^2+A_5^2$,\\
4. $A_7^2=c_{s2}(A_6^2+1)A_3^2+A_5^2$,\\
5. $A_8^2=A_3^2+A_7^2$.\\
Then
$$
\|\bar v_{sx}(kT)\|_{L_2}^2\le A_5^2
\leqno(3.5)
$$
and
$$
\|\bar v_{sx}(t)\|_{L_2}^2+c_{s1}\intop_{kT}^t\|\bar v_s(t')\|_{H^2}^2dt'\le 
A_8^2
\leqno(3.6)
$$
for all $t\in(kT,(k+1)T]$, $k\in\N_0$.

\Proof 
Differentiating $(1.5)_1$ with respect to $x$, multiplying by $\bar v_{sx}$ 
and integrating over $\Omega$ yields
$$
{1\over2}{d\over dt}\|\bar v_{sx}\|_{L_2}^2+\nu\|\bar v_{sxx}\|_{L_2}^2\le
\|\bar v_{sx}\|_{L_3}^3+\|\bar f_s\|_{L_2}\|\bar v_{sxx}\|_{L_2}.
$$
Using the Young inequality we get
$$
{1\over2}{d\over dt}\|\bar v_{sx}\|_{L_2}^2+{\nu\over2}
\|\bar v_{sxx}\|_{L_2}^2\le\|\bar v_{sx}\|_{L_3}^3+{1\over2\nu}
\|\bar f_s\|_{L_2}^2.
\leqno(3.7)
$$
Applying the interpolation inequality (see \cite{4})
$$
\|u\|_{L_3}\le c\|u_x\|_{L_2}^{1/3}\|u\|_{L_2}^{2/3}
$$
to the first term on the r.h.s. of (3.7), which holds for $\bar v_{sx}$ such 
that $\intop_\Omega\bar v_{sx}dx=0$, gives
$$
{d\over dt}\|\bar v_{sx}\|_{L_2}^2+\nu\|\bar v_{sxx}\|_{L_2}^2\le c_{s2}
\|\bar v_{sx}\|_{L_2}^4+c_{s2}\|\bar f_s\|_{L_2}^2.
\leqno(3.8)
$$
In view of inequality (2.5) we have
$$
{d\over dt}\|\bar v_{sx}\|_{L_2}^2+c_{s1}\|\bar v_{sx}\|_{L_2}^2\le c_{s2}
\|\bar v_{sx}\|_{L_2}^4+c_{s2}\|\bar f_s\|_{L_2}^2.
\leqno(3.9)
$$
Considering inequality (3.9) for $t\in[kT,(k+1)T]$ implies
$$\eqal{
&{d\over dt}\bigg(\|\bar v_{sx}\|_{L_2}^2e^{c_{s1}t-c_{s2}\intop_{kT}^t
\|\bar v_{sx}(t')\|_{L_2}^2dt'}\bigg)
\le c_{s2}\|f_s\|_{L_2}^2e^{c_{s1}t-c_{s2}\intop_{kT}^t
\|\bar v_{sx}(t')\|_{L_2}^2dt'}.\cr}
\leqno(3.10)
$$
Integrating (3.10) with respect to time from $t=kT$ to $t\in(kT,(k+1)T]$ 
we obtain
$$\eqal{
\|\bar v_{sx}(t)\|_{L_2}^2&\le e^{c_{s2}\intop_{kT}^t\|v_{sx}(t')\|_{L_2}^2dt'}
\cdot\intop_{kT}^t\|f_s(t')\|_{L_2}^2dt'\cr
&\quad+e^{-c_{s1}(t-kT)+c_{s2}\intop_{kT}^t\|v_{sx}(t')\|_{L_2}^2dt'}
\|\bar v_{sx}(kT)\|_{L_2}^2.\cr}
\leqno(3.11)
$$
Setting $t=(k+1)T$ in (3.11) and using (3.2) yields
$$\eqal{
\|\bar v_{sx}((k+1)T)\|_{L_2}^2&\le e^{c_{s2}A_3^2}\intop_{kT}^{(k+1)T}
\|f_s(t')\|_{L_2}^2dt'\cr
&\quad+e^{-c_{s1}T+c_{s2}A_3^2}\|\bar v_{sx}(kT)\|_{L_2}^2.\cr}
\leqno(3.12)
$$
In view of assumption (3.4) and notation 1. of the lemma we can write (3.12) 
briefly as
$$
\|\bar v_{sx}((k+1)T)\|_{L_2}^2\le A_4^2+e^{-{c_{s1}\over2}T}
\|v_{sx}(kT)\|_{L_2}^2.
$$
Hence iteration implies (3.5):
$$\eqal{
\|\bar v_{sx}(kT)\|_{L_2}^2&\le{A_4^2\over1-e^{-{c_{s1}\over2}T}}+
e^{-{c_{s1}\over2}kT}\|\bar v_{sx}(0)\|_{L_2}^2\cr
&\le{A_4^2\over1-e^{-{c_{s1}\over2}T}}+\|\bar v_{sx}(0)\|_{L_2}^2=A_5^2,\cr}
$$
where notation 2. is used. Employing (3.5) in (3.11) gives
$$\eqal{
\|\bar v_{sx}(t)\|_{L_2}^2& \le e^{c_{s2}A_3^2}\intop_{kT}^{(k+1)T}
\|\bar f_s(t')\|_{L_2}^2dt'+e^{-c_{s1}T+c_{s2}A_3^2}A_5^2\cr
&\le c_{s1}e^{c_{s2}A_3^2}A_1^2+A_5^2=A_4^2+A_5^2\equiv A_6^2\cr}
$$
for $t\in[kT,(k+1)T]$, where we used assumption 1. of Lemma 3.1 together with 
assumption (3.4) and notation 1. of the present lemma.

Integrating (3.8) with respect to time from $t=kT$ to $t\in(kT,(k+1)T]$ 
we obtain
$$\eqal{
&\|\bar v_{sx}(t)\|_{L_2}^2+\nu\intop_{kT}^t\|\bar v_{sxx}(t')\|_{L_2}^2dt'\le
c_{s2}\sup_t\|\bar v_{sx}(t)\|_{L_2}^2dt'\cr
&\le c_{s2}\sup_t\|\bar v_{sx}(t)\|_{L_2}^2\intop_{kT}^t
\|\bar v_{sx}(t')\|_{L_2}^2dt'+c_{s3}\intop_{kT}^t\|\bar f_s(t')\|_{L_2}^2dt'+
\|\bar v_{sx}(kT)\|_{L_2}^2\cr
&\le c_{s2}[A_6^2A_3^2+A_3^2]+A_5^2\equiv A_7^2.\cr}
$$
This implies (3.6) and ends the proof.
\kwadrat
\goodbreak

Inequalities (3.2) and (3.6) imply
$$
\|\bar v_s(t)\|_{H^1}^2+\intop_{kT}^t\|\bar v_s(t')\|_{H^2}^2dt'\le 
A_3^2+A_7^2\equiv A_8^2
\leqno(3.13)
$$
for all $t\in(kT,(k+1)T]$, $k\in\N_0$.

\proclaim Lemma 3.3. 
Suppose there exists a constant $A_9$ such that
$$
\sup_k\sup_{kT\le t\le(k+1)T}\big|\intop_0^t\diagintop_\Omega f_s(t')
dxdt'+\diagintop_\Omega v_s(0)dx\big|\le A_9<\infty.
$$
Let the assumptions of Lemmas 3.1 and 3.2 hold. Then there exists a solution 
to problem (1.5) such that $\bar v_s\in H^{2,1}(\Omega\times(kT,(k+1)T))$, 
$\nabla\bar p_s\in L_2(\Omega\times(kT,(k+1)T))$, $k\in\N_0$ and
$$\eqal{
&\|\bar v_s\|_{H^{2,1}(\Omega\times(kT,(k+1)T))}^2+
\|\nabla\bar p_s\|_{L_2(\Omega\times(kT,(k+1)T))}^2
\le cA_8^2(1+A_8^2)+cA_8^2A_9^2.\cr}
\leqno(3.14)
$$

\Proof 
Multiplying $(1.5)_1$ by $\bar v_{st}$, integrating over $\Omega$ and with 
respect to time from $kT$ to $(k+1)T$ gives
$$\eqal{
&\|\bar v_{st}\|_{L_2(\Omega\times(kT,(k+1)T))}^2\le c
\|\bar f_s\|_{L_2(\Omega\times(kT,(k+1)T))}^2+c\intop_{kT}^{(k+1)T}
\intop_\Omega|v_s|^2|\bar v_{sx}|^2dxdt\cr
&\quad+\|\bar v_{sx}(kT)\|_{L_2(\Omega)}^2\le cA_8^2(1+A_8^2+A_9^2),\cr}
$$
where
$$\eqal{
&\intop_{kT}^{(k+1)T}\intop_\Omega|v_s|^2|\bar v_{sx}|^2dxdt\le c
\|v_s\|_{L_\infty(kT,(k+1)T;H^1(\Omega))}^2
\|\bar v_s\|_{L_2(kT,(k+1)T;H^2(\Omega))}^2\cr
&\le c\bigg(\|\bar v_s\|_{L_\infty(kT,(k+1)T;H^1(\Omega))}^2+
\bigg\|\diagintop_\Omega v_sdx\bigg\|_{L_\infty(kT,(k+1)T;H^1(\Omega))}^2\bigg)
\|\bar v_s\|_{L_2(kT,(k+1)T;H^2(\Omega))}^2\cr
&\le c(A_8^2+A_9^2)A_8^2.\cr}
$$
Next, $(1.5)_1$ yields
$$\eqal{
&\|\nabla\bar p_s\|_{L_2(\Omega\times(kT,(k+1)T))}^2\le
\|\bar v_{st}\|_{L_2(\Omega\times(kT,(k+1)T))}^2\cr
&\quad+\|\bar v_s\|_{L_2(kT,(k+1)T;H^2(\Omega))}^2\cr
&\quad+c\|v_s\|_{L_\infty(kT,(k+1)T;H^1(\Omega))}^2
\|\bar v_s\|_{L_2(kT,(k+1)T;H^2(\Omega))}^2+
\|\bar f_s\|_{L_2(kT,(k+1)T;L_2(\Omega))}^2\cr
&\le cA_8^2+c(A_8^2+A_9^2)A_8^2.\cr}
$$
Hence (3.14) holds. Having estimate (3.14) existence follows by the 
Faedo-Galerkin method. This concludes the proof.
\kwadrat

\noindent
To prove stability of 2d solutions we need more regular 2d solutions than the 
one given in Lemma 3.2. Namely, we need

\proclaim Lemma 3.4. 
Let the assumptions of Lemma 3.2 be satisfied. Suppose that:\\
1. $A_{10}^2=\sup_k\intop_{kT}^{(k+1)T}\|\bar f_{sx}(t)\|_{L_2}^2dt<\infty$,\\
2. $A_{11}^2=\exp(c_{s3}A_8^2)c_{s3}A_{10}^2$, $c_{s4}={c_{s3}\over c_{s1}}$,\\
3. $T$ is so large that $-c_{s1}T/2+c_{s4}A_8^2\le0$,\\
4. $T$ is so large that $1-e^{c_{s1}T/2}\ge1/2$,\\
5. $A_{12}^2=2A_{11}^2+\|\bar v_{sxx}(0)\|_{L_2}^2<\infty$,\\
6. $A_{13}^2=A_{11}^2+A_{12}^2\exp(c_{s4}A_8^2)$,\\
7. $A_{14}^2=c_{s3}(A_{13}^2A_8^2+A_{10}^2)+A_{12}^2$,\\
where $c_{s3}>0$ is the constant from (3.19) below. Then
$$
\|\bar v_{sxx}(t)\|_{L_2}^2\le A_{13}^2,
\leqno(3.15)
$$
$$
\|\bar v_{sxx}(t)\|_{L_2}^2+c_{s1}\intop_{kT}^t
\|\bar v_{sxx}(t')\|_{H^1}^2dt'\le A_{14}^2
\leqno(3.16)
$$
for all $t\in(kT,(k+1)T]$, $k\in\N_0$.

\Proof 
Differentiating $(1.5)_1$ twice with respect to $x$, multiplying the result by 
$\bar v_{sxx}$, integrating over $\Omega$ and by parts yield
$$\eqal{
&{1\over2}{d\over dt}\|\bar v_{sxx}\|_{L_2}^2+\nu
\|\nabla\bar v_{sxx}\|_{L_2}^2=-\intop_\Omega\bar v_{sxx}\cdot\nabla\bar v_s
\cdot\bar v_{sxx}dx\cr
&\quad-2\intop_\Omega\bar v_{sx}\cdot\nabla\bar v_{sx}\cdot\bar v_{sxx}dx+
\intop_\Omega\bar f_{sxx}\cdot\bar v_{sxx}dx.\cr}
\leqno(3.17)
$$
Using the fact that $\bar v_s$ is divergence free we integrate by parts in the 
first two integrals on the r.h.s. of (3.17). We also integrate by parts in the 
third integral. Applying the H\"older and Young inequalities we obtain
$$\eqal{
&{1\over2}{d\over dt}\|\bar v_{sxx}\|_{L_2}^2+\nu
\|\nabla\bar v_{sxx}\|_{L_2}^2\le\varepsilon
\|\nabla\bar v_{sxx}\|_{L_2}^2\cr
&\quad+c(1/\varepsilon)\bigg(\intop_\Omega|\bar v_{sxx}|^2|\bar v_s|^2dx+
\intop_\Omega|\bar v_{sx}|^4dx+\intop_\Omega|\bar f_{sx}|^2dx\bigg).\cr}
$$
Hence for sufficiently small $\varepsilon$, from inequality (2.5) we get
$$
{d\over dt}\|\bar v_{sxx}\|_{L_2}^2+c_{s1}\|\bar v_{sxx}\|_{H^1}^2\le c
(\|\bar v_s\|_{L_\infty}^2\|\bar v_{sxx}\|_{L_2}^2
+\|\bar v_{sx}\|_{L_4}^4+\|\bar f_{sx}\|_{L_2}^2).
\leqno(3.18)
$$
Now, (3.18) implies that for $t\in(kT,(k+1)T)$
$$\eqal{
&{d\over dt}\bigg(\|\bar v_{sxx}\|_{L_2}^2\exp\bigg(c_{s1}t-c_{s3}\intop_{kT}^t
\|\bar v_s\|_{H^2}^2dt'\bigg)\bigg)\cr
&\le c_{s3}\|\bar f_{sx}\|_{L_2}^2\exp\bigg(c_{s1}t-
c_{s3}\intop_{kT}^t\|\bar v_s\|_{H^2}^2dt'\bigg),\cr}
\leqno(3.19)
$$
where we have used the estimates 
$\|\bar v_{sx}\|_{L_4}\le c\|\bar v_s\|_{H^2}$ and 
$\|\bar v_{sx}\|_{L_4}\le c\|\bar v_{sxx}\|_{L_2}$.

Integrating (3.19) with respect to time from $kT$ to $t\in(kT,(k+1)T]$, 
$k\in\N_0$, yields
$$\eqal{
&\|\bar v_{sxx}(t)\|_{L_2}^2\le c_{s3}\exp\bigg(-c_{s1}t+c_{s3}\intop_{kT}^t
\|\bar v_s(t')\|_{H^2}^2dt'\bigg)\cr
&\quad\cdot\intop_{kT}^t\|\bar f_{sx}(t')\|_{L_2}^2\exp
\bigg(c_{s1}t'-c_{s3}\intop_{kT}^{t'}\|\bar v_s(t'')\|_{H^2}^2dt''\bigg)dt'\cr
&\quad+\|\bar v_{sxx}(kT)\|_{L_2}^2\exp\bigg(-c_{s1}(t-kT)+c_{s3}\intop_{kT}^t
\|\bar v_s(t')\|_{H^2}^2dt'\bigg).\cr}
\leqno(3.20)
$$
Using notation 1. and (3.6) we obtain from (3.20) the inequality
$$\eqal{
&\|\bar v_{sxx}(t)\|_{L_2}^2\le c_{s3}\exp\bigg({c_{s3}\over c_{s1}}A_8^2\bigg)
A_{10}^2\cr
&\quad+\|\bar v_{sxx}(kT)\|_{L_2}^2\exp\bigg(-c_{s1}(t-kT)+{c_{s3}\over c_{s1}}
A_8^2\bigg).\cr}
$$
In view of notation 2. we have
$$
\|\bar v_{sxx}(t)\|_{L_2}^2\le A_{11}^2+\|\bar v_{sxx}(kT)\|_{L_2}^2\exp
(-c_{s1}(t-kT)+c_{s4}A_8^2).
\leqno(3.21)
$$
For $t=(k+1)T$, inequality (3.21) takes the form
$$
\|\bar v_{sxx}((k+1)T)\|_{L_2}^2\le A_{11}^2+\|\bar v_{sxx}(kT)\|_{L_2}^2\exp
(-c_{s1}T+c_{s4}A_8^2).
$$
Assumption 3. implies
$$
\|\bar v_{sxx}((k+1)T)\|_{L_2}^2\le A_{11}^2+e^{-c_{s1}T/2}
\|\bar v_{sxx}(kT)\|_{L_2}^2.
$$
Hence, iteration yields
$$\eqal{
\|\bar v_{sxx}(kT)\|_{L_2}^2&\le{A_{11}^2\over1-e^{-c_{s1}T/2}}+
e^{-c_{s1}kT/2}\|\bar v_{sxx}(0)\|_{L_2}^2\cr
&\le2A_{11}^2+\|\bar v_{sxx}(0)\|_{L_2}^2=A_{12}^2,\cr}
\leqno(3.22)
$$
where Assumption 4 is utilized. Employing (3.22) in (3.21) gives (3.15).

Integrating (3.18) with respect to time implies the estimate
$$
\|\bar v_{sxx}(t)\|_{L_2}^2+c_{s1}\intop_{kT}^t
\|\bar v_{sxx}(t')\|_{H^1}^2dt'\le c_{s3}A_{13}^2A_8^2+c_{s3}A_{10}^2+
A_{12}^2=A_{14}^2
$$
for all $t\in(kT,(k+1)T]$, $k\in\N_0$. This implies (3.16) and concludes the 
proof.
\kwadrat
\goodbreak

\Remark{3.5.} 
Applying Faedo-Galerkin approximations and using Lemmas 3.1, 3.2, 3.4 and 
estimates (3.15)--(3.16), we conclude that the assertion of Theorem 1.1 
holds. Employing additionally Lemma 3.3 we obtain Theorem 1.2.

\section{4. Stability}

To prove the stability of two-dimensional solutions we have to find solutions 
to problem (1.3) such that the inequality $\|u(0)\|_{H^1}\le\gamma$ 
implies that $\|u(t)\|_{H^1}\le c\gamma$ for $\gamma$ sufficiently small and 
for all $t\in\R_+$, where $c>0$ is a constant.

First we derive an energy type estimate for solutions to problem (1.6).

\proclaim Lemma 4.1. 
Let the assumptions of Lemmas 3.1, 3.2 hold. 
Assume that \\
$g\in L_2(kT,(k+1)T;L_2(\Omega))$, $k\in\N_0$ and $\bar u$ 
satisfies (1.6). Assume that\\
1. $B_1^2=\sup_k\intop_{kT}^{(k+1)T}\|\bar g(t')\|_{L_2}^2dt'<\infty$,\\
2. $B_2^2=\sup_k\intop_{kT}^{(k+1)T}\big|\intop_0^t\diagintop_\Omega gdxdt'+
\diagintop_\Omega u(0)dx\big|^2dt<\infty$,\\
3. 
$B_3^2=\big(c_2B_1^2+c_2A_3^2B_2^2\big)\exp(c_2A_8^2)$, where $c_1>0$ is the 
constant from (2.6) and $c_2>0$ appears in (4.3).\\
4. $T$ is so large that $-c_1T/2+A_8^2\le0$, $A_8^2$ appears in (3.13),\\
5. $T$ is so large that $1-e^{-c_1T/2}\ge1/2$,\\
6. $B_4^2=B_3^2+\exp(c_2A_8^2)(2B_4^2+\|\bar u(0)\|_{L_2}^2)$.\\
Then
$$
\|\bar u(t)\|_{L_2}^2+c_1\intop_{kT}^t\|\bar u(t')\|_{H^1}^2dt'\le
c_2A_8^2B_4^2+c_2A_3^2B_2^2+c_2B_1^2+B_3^2\equiv B_5^2
\leqno(4.1)
$$
for all $t\in(kT,(k+1)T)$, $k\in\N_0$.

\Proof 
Multiplying $(1.6)_1$ by $\bar u$, integrating over $\Omega$, by parts and 
using the periodic boundary conditions we obtain
$$\eqal{
&{1\over2}{d\over dt}\|\bar u\|_{L_2}^2+\nu\|\nabla\bar u\|_{L_2}^2=
-\intop_\Omega u\cdot\nabla\bar v_s\cdot\bar udx+\intop_\Omega
\bar g\cdot\bar udx\cr
&=-\intop_\Omega(\bar u+\diagintop_\Omega udx)\cdot\nabla\bar v_s\cdot\bar udx+
\intop_\Omega\bar g\cdot\bar udx\cr
&=-\intop_\Omega\bar u\cdot \nabla\bar v_s\cdot\bar udx+\diagintop_\Omega
udx\cdot\intop_\Omega\bar v_s\cdot\nabla\bar udx+\intop_\Omega\bar g\cdot
\bar udx.\cr}
\leqno(4.2)
$$
Using the H\"older and Young inequalities we get
$$\eqal{
&{1\over2}{d\over dt}\|\bar u\|_{L_2}^2+\nu\|\nabla\bar u\|_{L_2}^2\le
\varepsilon(\|\bar u\|_{L_6}^2+\|\bar u_x\|_{L_2}^2+\|\bar u\|_{L_2}^2)\cr
&\quad+c(\varepsilon)\bigg(\|\bar v_{sx}\|_{L_3}^2\|\bar u\|_{L_2}^2+
\|\bar v_s\|_{L_2}^2\bigg|\intop_0^t\diagintop_\Omega gdxdt'+
\diagintop_\Omega u(0)dx\bigg|^2+\|\bar g\|_{L_2}^2\bigg).\cr}
$$
Assuming that $\varepsilon$ is sufficiently small and applying 
inequality (2.6) yields
$$\eqal{
&{d\over dt}\|\bar u\|_{L_2}^2+c_1\|\bar u\|_{H^1}^2\le c_2
\|\bar v_{sx}\|_{L_3}^2\|\bar u\|_{L_2}^2\cr
&\quad+c_2\|\bar v_s\|_{L_2}^2\bigg|\intop_0^t\diagintop_\Omega gdxdt'+
\diagintop_\Omega u(0)dx\bigg|^2+c_2\|\bar g\|_{L_2}^2,\cr}
\leqno(4.3)
$$
where $c_1$ is the constant from (2.6). Inequality (4.3) implies
$$\eqal{
&{d\over dt}\big(\|\bar u\|_{L_2}^2
e^{c_1t-c_2\intop_{kT}^t\|\bar v_s\|_{H^2}^2dt'}\big)\le c_2
\bigg(\|\bar g\|_{L_2}^2\cr
&\quad+\|\bar v_s\|_{L_2}^2\bigg|\intop_0^t\diagintop_\Omega gdxdt'+
\diagintop_\Omega u(0)dx\bigg|^2\bigg)
e^{c_1t-c_2\intop_{kT}^t\|\bar v_s(t')\|_{H^2}^2dt'}\cr}
\leqno(4.4)
$$
for all $t\in(kT,(k+1)T]$, $k\in\N_0$.

Integrating (4.4) with respect to time from $kT$ to $t\in(kT,(k+1)T]$ yields
$$\eqal{
&\|\bar u(t)\|_{L_2}^2\le c_2\exp\bigg(-c_1t+c_2\intop_{kT}^t
\|\bar v_s(t')\|_{H^2}^2dt'\bigg)\cdot\cr
&\quad\cdot\intop_{kT}^t\bigg(\|\bar g(t')\|_{L_2}^2+A_3^2\bigg|\intop_0^{t'}
\diagintop_\Omega gdxdt''+\diagintop_\Omega u(0)dx\bigg|^2\bigg)e^{c_1t'}dt'\cr
&\quad+\exp\bigg(-c_1(t-kT)+c_2\intop_{kT}^t\|\bar v_s(t')\|_{H^2}^2dt'\bigg)
\|\bar u(kT)\|_{L_2}^2,\cr}
\leqno(4.5)
$$
where (3.2) is used. In view of Assumptions 1--4 and (3.13) we have
$$
\|\bar u(t)\|_{L_2}^2\le B_3^2+e^{-c_1(t-kT)+c_2A_8^2}\|\bar u(kT)\|_{L_2}^2.
$$
Setting $t=(k+1)T$ and using Assumption 5 we get
$$
\|\bar u((k+1)T)\|_{L_2}^2\le B_3^2+e^{-c_1T/2}\|\bar u(kT)\|_{L_2}^2.
$$
Hence, iteration implies
$$
\|\bar u(kT)\|_{L_2}^2\le{B_3^2\over1-e^{-c_1T/2}}+e^{-c_1kT/2}
\|\bar u(0)\|_{L_2}^2
\le2B_3^2+\|\bar u(0)\|_{L_2}^2,
\leqno(4.6)
$$
where Assumption 6 is used. Inserting (4.6) in (4.5) yields
$$\eqal{
\|\bar u(t)\|_{L_2}^2&\le B_3^2+e^{-c_1(t-kT)+c_2A_8^2}(2B_3^2+
\|\bar u(0)\|_{L_2}^2)\cr
&\le B_3^2+e^{c_2A_8^2}(2B_4^2+\|\bar u(0)\|_{L_2}^2)\equiv B_4^2.\cr}
$$
Integrating (4.3) with respect to time from $kT$ to $t\in(kT,(k+1)T]$ we derive
$$\eqal{
&\|\bar u(t)\|_{L_2}^2+c_1\intop_{kT}^t\|\bar u(t')\|_{H^1}^2dt'\le 
c_2A_8^2B_4^2+c_2A_3^2B_2^2+c_2B_1^2+B_3^2.\cr}
$$
This implies (4.1) and concludes the proof.
\kwadrat

Now, we show that the 3d solution to (1.1) remains close to the 2d solution of 
(1.2) if they are sufficiently close at the initial time.

\proclaim Lemma 4.2. 
Let the assumptions of Lemma 4.1 hold. Let $\gamma_*$ be so small that 
$c_1-c_3\gamma_*^4\ge c_1/2$, where $c_1>0$ is the constant from (2.6) and 
$c_3>0$ occurs in (4.10)--(4.11). Let $\gamma\in(0,\gamma_*]$. Assume that
$$
\|\bar u(0)\|_{H^1}^2\le\gamma,
$$
$$
c_3\|v_{sx}\|_{L_3}^2\bigg[\|\bar v_{sx}(t)\|_{L_3}^2B_5^2+\bigg|\intop_0^t
\diagintop_\Omega gdxdt+\diagintop_\Omega u(0)dx\bigg|^2\bigg]+c_3
\|\bar g\|_{L_2}^2\le{c_1\over4}\gamma \quad {\sl for\ all}\ \ t\in\R_+.
$$
Then
$$
\|\bar u(t)\|_{H^1}^2\le\gamma\quad {\sl for\ any}\ \ t\in\R_+.
$$

\noindent
This means that the 3d solution to (1.1) remains close to the 2d solution of 
(1.2) if their initial data and the external forces  for all time are 
sufficiently close.

\Proof 
Differentiating $(1.6)_1$ with respect to $x$, multiplying the result by 
$\bar u_x$, integrating over $\Omega$, by parts and employing the periodic 
boundary conditions we obtain
$$\eqal{
&{1\over2}{d\over dt}\|\bar u_x\|_{L_2}^2+\nu\|\bar u_{xx}\|_{L_2}^2\le
\|\bar u_x\|_{L_3}^3+\bigg|\intop_\Omega\bar v_{sx}\cdot\nabla u\cdot
\bar u_xdx\bigg|\cr
&\quad+2\bigg|\intop_\Omega\bar u_x\cdot\nabla\bar v_s\cdot\bar u_xdx\bigg|+
\bigg|\intop_\Omega u\cdot\nabla\bar v_s\cdot\bar u_{xx}dx\bigg|+
\bigg|\intop_\Omega\bar g\cdot\bar u_{xx}dx\bigg|.\cr}
\leqno(4.7)
$$
Adding (4.2) and (4.7), applying the H\"older and Young inequalities, we derive
$$\eqal{
&{d\over dt}\|\bar u\|_{H^1}^2+c_1\|\bar u\|_{H^2}^2\le c
(\|\bar u_x\|_{L_3}^3+\|\bar v_{sx}\|_{L_3}^2\|\bar u_x\|_{L_2}^2\cr
&\quad+\|\bar v_{sx}\|_{L_3}^2\|\bar u_x\|_{L_2}^2+
\|u\|_{L_6}^2\|\bar v_{sx}\|_{L_3}^2+\|\bar g\|_{L_2}^2).\cr}
$$
Using $\|u\|_{L_6}^2\le\|\bar u\|_{L_6}^2+\big|\diagintop_\Omega udx\big|^2$ 
and $\|\bar u\|_{L_6}\le c\|\bar u\|_{H^1}\le c\|\bar u_x\|_{L_2}$, which 
holds in view of the Poincar\'e inequality, we get
$$\eqal{
&{d\over dt}\|\bar u\|_{H^1}^2+c_1\|\bar u\|_{H^2}^2\le c
\bigg[\|\bar u_x\|_{L_3}^3+\|\bar v_{sx}\|_{L_3}^2\|\bar u_x\|_{L_2}^2
+\|v_{sx}\|_{L_3}^2\bigg|\diagintop_\Omega udx\bigg|^2+
\|\bar g\|_{L_2}^2\bigg].\cr}
\leqno(4.8)
$$
In view of (2.2) and the interpolation inequality (see [4, Ch. 3, Sect. 15])
$$
\|\bar u_x\|_{L_3}\le c\|\bar u_{xx}\|_{L_2}^{1/2}\|\bar u_x\|_{L_2}^{1/2}
$$
(which holds without the lower order term because 
$\intop_\Omega\bar u_xdx=0$), we obtain from (4.8) the inequality
$$\eqal{
&{d\over dt}\|\bar u\|_{H^1}^2+c_1\|\bar u\|_{H^2}^2\le c
\|\bar u_x\|_{L_2}^6+c\|\bar v_{sx}\|_{L_3}^2\|\bar u_x\|_{L_2}^2\cr
&\quad+c\|v_{sx}\|_{L_3}^2\bigg(\bigg|\intop_0^t\diagintop_\Omega 
gdxdt'+\diagintop_\Omega u(0)dx\bigg|^2\bigg)+c\|\bar g\|_{L_2}^2.\cr}
\leqno(4.9)
$$
Employing the interpolation inequality (see \cite[Ch. 3, Sect. 10]{4})
$$
\|\bar u_x\|_{L_2}\le\varepsilon^{1/2}\|\bar u_{xx}\|_{L_2}+
c\varepsilon^{-1/2}\|\bar u\|_{L_2}
$$
in (4.9) implies
$$\eqal{
&{d\over dt}\|\bar u\|_{H^1}^2+c_1\|\bar u\|_{H^2}^2\le c_3
\|\bar u_x\|_{L_2}^6\cr
&\quad+c_3\|\bar v_{sx}\|_{L_3}^2\|\bar u\|_{L_2}^2+c_3
\|\bar v_{sx}\|_{L_3}^2\bigg(\bigg|\intop_0^t\diagintop_\Omega 
gdxdt'+\diagintop_\Omega u(0)dx\bigg|^2\bigg)+c_3\|\bar g\|_{L_2}^2.\cr}
\leqno(4.10)
$$
In view of (4.1) we have $\|\bar u(t)\|_{L_2}\le B_6$.
Hence we can introduce the quantities:
$$\eqal{
G^2(t)&=c_3\|\bar v_{sx}(t)\|_{L_3}^2\bigg[B_6^2+
\bigg|\intop_0^t\diagintop_\Omega gdxdt'+\diagintop_\Omega u(0)dx\bigg|^2
\bigg]+c_4\|\bar g(t)\|_{L_2}^2,\cr
X(t)&=\|\bar u(t)\|_{H^1},\quad Y(t)=\|\bar u(t)\|_{H^2}.\cr}
$$
Then (4.10) takes the form
$$
{d\over dt}X^2\le-c_1Y^2+c_3X^4X^2+G^2.
$$
Since $X\le Y$ we have
$$
{d\over dt}X^2\le-X^2(c_1-c_3X^4)+G^2.
\leqno(4.11)
$$
Let $\gamma\in(0,\gamma_*]$, where $\gamma_*$ is so small that
$$
c_1-c_3\gamma_*^4\ge c_1/2.
\leqno(4.12)
$$
By the assumptions of the lemma,
$$
X^2(0)\le\gamma,\quad G^2(t)\le c_1{\gamma\over4}\quad {\rm for\ all}\ \ 
t\in\R_+.
$$
Suppose that
$$
t_*=\inf\{t\in\R_+:\ X^2(t)>\gamma\}>0.
$$
Then by (4.12) for $t\in(0,t_*]$ inequality (4.11) takes the form
$$
{d\over dt}X^2\le-{c_1\over2}X^2+G^2(t).
\leqno(4.13)
$$
Clearly, we have
$$
X^2(t_*)=\gamma\quad {\rm and}\ X^2(t)>\gamma\quad {\rm for}\ t>t_*.
\leqno(4.14)
$$
Then (4.13) yields
$$
{d\over dt}X^2(t)\bigg|_{t=t_*}\le c_1\bigg(-{\gamma\over2}+
{\gamma\over4}\bigg)<0
$$
contradicting with (4.14). Therefore
$$
X^2(t)<\gamma\quad {\rm for}\ \ t\in\R_+.
\leqno(4.15)
$$
This concludes the proof.
\kwadrat

\proclaim Lemma 4.3. 
Let $A_8$ be as introduced in (3.13), $A_9$ as in Lemma 3.3 and $\gamma$ as in 
Lemma 4.2. Let $T$ be as defined in Lemma 4.1. Let
$$
B_6=\sup_k\sup_{kT\le t\le(k+1)T}\bigg|\intop_0^t\diagintop_\Omega
g(t)dxdt+\diagintop_\Omega u(0)dx\bigg|.
$$
Then there exists a solution to problem (1.6) such that 
$\bar u\in H^{2,1}(\Omega\times(kT,(k+1)T))$, 
$\nabla\bar q\in L_2(kT,(k+1)T;L_2(\Omega))$, $k\in\N_0$, and
$$\eqal{
&\|\bar u\|_{H^{2,1}(\Omega\times(kT,(k+1)T))}^2+
\|\nabla\bar q\|_{L_2(kT,(k+1)T;L_2(\Omega))}^2\cr
&\le c[(T+1)\gamma^2+B_6^2][A_8^2(1+A_8^2+A_9^2)+A_9^2+(T+1)\gamma^2]+
c\gamma^2\equiv B_7^2.\cr}
$$

\Proof 
In view of the definition of $G$ we express (4.10) in the form
$$
{d\over dt}\|\bar u\|_{H^1}^2+c_1\|\bar u\|_{H^2}^2\le c_3
\|\bar u_x\|_{L_2}^6+G^2.
\leqno(4.16)
$$
Integrating (4.16) with respect to time from $kT$ to $t\in[kT,(k+1)T]$ and 
using (4.15) we derive
$$\eqal{
&\|\bar u(t)\|_{H^1}^2+c_1\intop_{kT}^t\|\bar u(t')\|_{H^2}^2dt'\le
c_3\gamma^6T+\gamma^2T+\gamma^2\le c(T+1)\gamma^2,\cr}
\leqno(4.17)
$$
because $\gamma<1$. Multiplying (1.6) by $\bar u_t$ and integrating over 
$\Omega\times(kT,(k+1)T)$ yields
$$\eqal{
&\intop_{kT}^{(k+1)T}\|\bar u_t(t)\|_{L_2}^2dt+\nu
\|\nabla\bar u((k+1)T)\|_{L_2}^2\cr
&\le c\intop_{kT}^{(k+1)T}(\|u\cdot\nabla\bar u\|_{L_2}^2+
\|v_s\cdot\nabla\bar u\|_{L_2}^2+\|u\cdot\nabla\bar v_s\|_{L_2}^2)dt\cr
&\quad+c\intop_{kT}^{(k+1)T}\|\bar g(t)\|_{L_2}^2dt+\nu
\|\nabla\bar u(kT)\|_{L_2}^2.\cr}
\leqno(4.18)
$$
The first term on the r.h.s. of (4.18) is estimated by
$$\eqal{
&c\sup_{kT\le t\le(k+1)T}\|u(t)\|_{H^1}^2(\|\bar u\|_{L_2(kT,(k+1)T;H^2)}^2+
\|\bar v_s\|_{L_2(kT,(k+1)T;H^2)}^2)\cr
&\quad+c\sup_{kT\le t\le(k+1)T}\|v_s(t)\|_{H^1}^2
\|\bar u\|_{L_2(kT,(k+1)T;H^2)}^2\cr
&\le c(\gamma^2+B_6^2)[A_8^2(1+A_8^2)+A_8^2A_9^2+(T+1)\gamma^2]+
c(A_8^2+A_9^2)(T+1)\gamma^2\cr
&\le B_7^2.\cr}
$$
Similarly,
$$\eqal{
&\|\nabla\bar q\|_{L_2(kT,(k+1)T;L_2)}\le c
\|\bar u\|_{H^{2,1}(\Omega\times(kT,(k+1)T))}\cr
&\quad+c\|u\cdot\nabla\bar u\|_{L_2(kT,(k+1)T;L_2)}+
c\|v_s\cdot\nabla\bar u\|_{L_2(kT,(k+1)T,L_2)}\cr
&\quad+c\|u\cdot\nabla\bar v_s\|_{L_2(kT,(k+1)T;L_2)}+
\|\bar g\|_{L_2(kT,(k+1)T;L_2)}\le cB_7.\cr}
$$
This concludes the proof.
\kwadrat
\goodbreak 

Now, we can complete the proofs of Theorems 1.3 and 1.4.
\vskip6pt

\noindent
{\bf The proofs of Theorems 1.3 and 1.4}

Inequalities (1.7) and (1.9) follow from Lemmas 4.2 and 4.3, respectively. 
The existence of solution is a consequence of applying the 
Faedo-Galerkin method and inequalities (1.7)--(1.8). Thus, we get the 
assertion of Theorem 1.3. Theorem 1.4 follows  from Lemma 4.3.
\kwadrat

\section{References}

\item{[1]} Auscher P., Dubois S. and Tchamitchian P.: On the stability of 
global solutions to Navier-Stokes equations in the space, Journal de 
Math\'ematiques Pures et Appliques 83 (2004), 673--697.

\item{[2]} Beir\~ao da Veiga H. and Secchi P.: $Lp$-stability for the strong 
solutions of the Navier-Stokes equations in the whole space, Arch. Ration. 
Mech. Anal. 98 (1987), 65--69.

\item{[3]} Bardos C., Lopes Filho M. C., Niu D., Nussenzveig Lopes H. J. and 
Titi E. S.: Stability of two-dimensional viscous incompressible flows under 
three-dimesional perturbations and inviscid symmetry breaking, SIAM J. Math. 
Anal. 45 (2013), 1871--1885.

\item{[4]} Besov, O. V.; Il'in, V. P.; Nikol'skii, S. M.: Integral 
representations of functions and imbedding theorems, Nauka, Moscow 1975 
(in Russian).

\item{[5]} Chemin J. I, and Gallagher I.: Wellposedness and stability results 
for the Navier-Stokes equations in $R^3$, Ann. I. H. Poincar\'e, Analyse Non 
Lin\'eaire, 26 (2009), 599--624.

\item{[6]} Gallagher I.: The tridimensional Navier-Stokes equations with 
almost bidimensional data: stability, uniqueness and life span, Internat. Mat. 
Res. Notices 18 (1997), 919--935.

\item{[7]} Gallagher I., Iftimie D. and Planchon F.: Asymptotics and stability 
for global solutions to the Navier-Stokes equations, Annales de I'Institut 
Fourier 53 (2003), 1387--1424.

\item{[8]} Iftimie D.: The 3d Navier-Stokes equations seen as a perturbation 
of the 2d Navier-Stokes equations, Bull. Soc. Math. France 127 (1999), 
473--517.

\item{[9]} Karch G. and Pilarczyk D.: Asymptotic stability of Landau 
solutions to Navier-Stokes system, Arch. Rational Mech. Anal. 202 (2011), 
115--131.

\item{[10]} Mucha P. B.: Stability of nontrivial solutions of the Navier-Stokes 
system on the three-dimensional torus, J. Differential Equations 172 (2001), 
359--375.

\item{[11]} Mucha P. B.: Stability of constant solutions to the Navier-Stokes 
system in $\R^3$, Appl. Math. 28 (2001), 301--310.

\item{[12]} Mucha P. B.: Stability of 2d incompressible flows in $R^3$, J. 
Diff. Eqs. 245 (2008), 2355--2367.

\item{[13]} Ponce G., Racke R., Sideris T. C. and Titi E. S.: Global stability 
of large solutions to the 3d Navier-Stokes equations, Comm. Math. Phys. 159 
(1994), 329--341.

\item{[14]} Zhou Y.: Asymptotic stability for the Navier-Stokes equations 
in the marginal class, Proc Roy. Soc. Edinburgh 136 (2006), 1099--1109.

\bye